\documentclass{jms}

\newcommand{\abs}[1]{\left|#1\right|}
\renewcommand{\d}{\delta}
\newcommand{\eps}{\varepsilon}
\renewcommand{\le}{\leqslant}
\renewcommand{\ge}{\geqslant}

\newcommand{\lyapd}{L_{2+\delta,n}}
\newcommand{\I}{{\bf1}}

\paperTitle{Lower bounds for the constants in non-uniform estimates of the rate of convergence in the CLT}

\authorsShort{I.\,G.~Shevtsova}
\authorsFull{I.\,G. Shevtsova\first\also\second\also\third}
\addAuthorInfo{Hangzhou Dianzi University, Hangzhou, China, e-mail: \url{ishevtsova@cs.msu.ru}}
\addAuthorInfo{Lomonosov Moscow State University, Moscow, Russia}
\addAuthorInfo{FRSC ``Informatics and Control''  of RAS, Moscow, Russia}

\paperAbstract{We conduct a comparative analysis of the constants in the Nagaev--Bikelis and Bikelis--Petrov inequalities which establish non-uniform estimates of the rate of convergence in the central limit theorem for sums of independent random variables possessing finite absolute moments of order $2+\d$ with $\d\in[0,1]$. We provide lower bounds for the above constants and also for the constants in the structural improvements of Nagaev--Bikelis' inequality. The lower bounds in Nagaev--Bikelis' inequality and it's structural improvements are given in dependence on $\d$ and a structural parameter $s$ as well as uniform with respect to both $\d$ and $s$. Lower bounds for the constants in Nagaev--Bikelis'  with $\d<1$ and Bikelis--Petrov's inequalities are presented for the first time.}

\begin{document}
\maketitle

\section{Introduction}

Let $X_1,X_2,\ldots,X_n$ be independent random variables (r.v.'s) with distribution functions (d.f.'s) $F_1,F_2,\ldots,F_n,$ and $\E X_k=0,$ $\sigma_k^2:=\E X_k<\infty$,
$$
S_n:=\sum_{k=1}^nX_k,\quad B_n^2:=\sum_{k=1}^n\sigma_k^2=\D S_n>0.
$$
Let us denote
$$
\overline F_n(x):= \Prob(S_n<xB_n),\quad \Phi(x)=\frac1{\sqrt{2\pi}}\int_{-\infty}^xe^{-t^2/2}dt,
\quad
\Delta_n(x)=|\overline F_n(x)-\Phi(x)|,\quad  x\in\R,
$$
$$
L_n(\eps):= \frac1{B_n^2}\sum_{k=1}^n\E X_k^2\I(|X_k|>\eps B_n), \quad \eps>0,\ n\in\N.
$$
If
$$
\max_{1\le k\le n}\E|X_{k}|^{2+\d}<\infty\quad \text{for some }\d>0
$$
denote also
$$
\lyapd:= \frac1{B_n^{2+\d}}\sum_{k=1}^n\E|X_{k}|^{2+\d},\quad \d\in(0,1].
$$
We shall also use the notation $\lyapd:=1$ for $\d=0$. The quantities $L_n(\,\cdot\,)$ and $\lyapd$ are called the Lindeberg fraction and the Lyapounov fraction, respectively.

In the case of independent and identically distributed (i.i.d.) random summands $X_1,\ldots,X_n$ and $\d=1$, Nagaev~\cite{Nagaev1965} proved that
\begin{equation}\label{NagaevBikelisIneq}
\Delta_n(x)\le K_0\frac{\lyapd}{1+|x|^{2+\d}},\quad x\in\R,\ n\in\N,
\end{equation}
where $K_0=K_0(\d)$ is an absolute constant for every value of $\d\in(0,1]$ and may be chosen even independent on $\d$, i.e. to be an absolute constant \textit{uniformly} for all $\d\in[0,1]$. Inequality~\eqref{NagaevBikelisIneq} was proved in the general situation (in the presented form) by  Bikelis~\cite{Bikelis1966} one year later, i.\,e., in 1966.

The first upper bounds for $K_0(\d)$ (in dependence on $\d$) were obtained by Paditz~\cite{Paditz1976,Paditz1977,Paditz1978,Paditz1979} in 1976--1979 and were of order $100\div2000$. These estimates further were considered and consequently improved by Michel~\cite{Michel1981} (i.i.d. case and $\d=1$) in 1981, Tysiak~\cite{Tysiak1983} in 1983, Mirachmedov and Paditz~\cite{Mirachmedov1984,PaditzMirachmedov1986,Paditz1989}  in 1984--1989, Nefedova, Shevtsova, Grigorieva and Popov~\cite{NefedovaShevtsova2011,NefedovaShevtsova2012, GrigorievaPopov2012SMI,Shevtsova2013Inf,Shevtsova2016,Shevtsova2017book} in 2011--2017. A detailed survey may be found, e.\,g., in~\cite{Shevtsova2017book}.  The best known upper bounds for  $K_0(\d)$ are obtained in~\cite{Shevtsova2017book} (see also  announcement~\cite{Shevtsova2013Inf}) and are presented in table\,\ref{TabBEnonuni1+x_struct} in the second and fifths columns.

\begin{table}[h]
\centering
\begin{tabular}{||c||c|c|c||c|c|c||}
\hline &\multicolumn{3}{c||}{non-i.i.d. case}
&\multicolumn{3}{c||}{i.i.d. case}\\
\cline{2-7} $\delta$ & \vphantom{$\frac{1^1}2$}
$K_0(\d)$ & $K_{s_1}(\d)$ & $s_1$ &
$K_0(\d)$ & $K_{s_1}(\d)$ & $s_1$
\\\hline
$1$  &21.82&18.19&1    &17.36&15.70&0.646\\
$0.9$&20.07&16.65&1    &16.24&14.61&0.619\\
$0.8$&18.53&15.34&1    &15.20&13.61&0.625\\
$0.7$&17.14&14.20&1    &14.13&12.71&0.570\\
$0.6$&15.91&13.19&0.859&13.15&11.90&0.498\\
$0.5$&14.84&12.30&0.834&12.26&11.17&0.428\\
$0.4$&13.92&11.53&0.806&11.43&10.51&0.350\\
$0.3$&13.10&10.86&0.778&10.66&9.93 &0.273\\
$0.2$&12.35&10.28&0.748&9.92 &9.42 &0.183\\
$0.1$&11.67& 9.77&0.710&9.18 &8.97 &0.074\\
\hline
\end{tabular}
\caption{Upper bounds for the constants $K_s(\delta)$ from inequalities~\eqref{NagaevBikelisIneq} and~\eqref{BEnonuni1+x_struct} for some $s\in[0,1]$ and $\delta\in(0,1]$.}
\label{TabBEnonuni1+x_struct}
\end{table}

Let us note that Bikelis in~\cite[Theorem 4]{Bikelis1966} obtained, in fact, a stronger result, which can be called a non-uniform analogue of the Osipov inequality~\cite{Osipov1966}: for all $x\in\R$ and  $n\in\N$ 
\begin{eqnarray}\label{BikelisIneqInt}
\Delta_n(x)&\le& \frac{A}{(1+|x|)^3B_n}\int_0^{(1+|x|)B_n}L_n(z) d z 
\\
\label{BikelisIneqMin}
&=&\frac{A} {(1+|x|)^3B_n^3}\sum_{k=1}^n\E X_k^2\min\big\{|X_k|,\,(1+|x|)B_n\big\}
\end{eqnarray}
\begin{equation}\label{BikelisIneq}
=A\sum_{k=1}^n\bigg[\frac{\E X_k^2\I(|X_k|>(1+|x|)B_n)}{(1+|x|)^2B_n^2}
+ \frac{\E|X_k|^3\I(|X_k|\le (1+|x|)B_n)}{(1+|x|)^3B_n^3}\bigg],
\end{equation}
where $A$ is an absolute constant. Inequality~\eqref{BikelisIneqInt} trivially yields~\eqref{NagaevBikelisIneq}: Indeed, following the reasoning of Bikelis and multiplying the quadratic and the cubic functions in the integrands in~\eqref{BikelisIneq} by
$$
\Big(\frac{|X_k|}{(1+|x|)B_n}\Big)^\d\ge1 \quad \text{and}\quad \Big(\frac{(1+|x|)B_n}{|X_k|}\Big)^{1-\d}\ge1,
$$
respectively, we get~\eqref{NagaevBikelisIneq} with
\begin{equation}\label{K0<=A*2^(1+d)}
K_0(\d)\le A\sup_{x>0}\frac{(1+x)^{2+\d}}{1+x^{2+\d}}= A\cdot\frac{(1+x)^{2+\d}}{1+x^{2+\d}}\bigg|_{x=1}= A\cdot 2^{1+\d}\le 4A,\quad \d\in[0,1].
\end{equation}
Though expression~\eqref{BikelisIneq}, obviously, was kept in mind by Bikelis when he was deducing~\eqref{NagaevBikelisIneq},  it was not given in the explicit form in~\cite{Bikelis1966}, however. Inequality~\eqref{BikelisIneq} appears for the first time only in the work of Petrov~\cite{Petrov1979} in 1979, where the author deduces it from Bikelis' inequality~\eqref{BikelisIneqInt} in the course of the proof of inequality~\eqref{PetrovIneq} below. 

In order to formulate Petrov's inequality~\eqref{PetrovIneq} below let us introduce a set $\mathcal G$ of all even functions $g\colon\R\to\R_+$ such that $g(x)>0$ for $x>0$ and the functions $g(x)$, $x/g(x)$ are non-decreasing for $x>0.$ In recent paper~\cite{GabdullinMakarenkoShevtsova2018} it was proved that every function $g\in\mathcal G$ satisfies the following inequalities for all $x\in\R\setminus\{0\}$ and $a>0:$
\begin{equation}\label{KatzPetrovGproperties}
g_*(x;a):= \min\Big\{1,\frac{|x|}a\Big\}\ \le\ \frac{g(x)}{g(a)}\ \le\ \max\Big\{1,\frac{|x|}a\Big\}=: g^*(x;a),
\end{equation}
where $g_*(\,\cdot\,;a),g^*(\,\cdot\,;a)\in\mathcal G$ for every $a>0$.

Petrov~\cite{Petrov1979} proved that for every function $g\in\mathcal G$ such that $\max\limits_{1\le k\le n}\E X_k^2g(X_k)<\infty$ we have
\begin{equation}\label{PetrovIneq}
\Delta_n(x)\le \frac{A\sum_{k=1}^n\E X_k^2g(X_k)}{(1+|x|)^2B_n^2g((1+|x|)B_n)},\quad x\in\R,\ n\in\N,
\end{equation}
with the same constant  $A$ as in~\eqref{BikelisIneq}, that is, \textit{universal} for all functions $g\in\mathcal G$.

Observe that inequality~\eqref{PetrovIneq}  follows trivially from~\eqref{BikelisIneq} with the account of~\eqref{KatzPetrovGproperties}. The proofs of~\eqref{PetrovIneq} in~\cite{Petrov1979} and~\eqref{KatzPetrovGproperties} in~\cite{GabdullinMakarenkoShevtsova2018} are based on the same ideas, Petrov applying them only to the concrete functional of a function $g\in\mathcal G$, while the authors of~\cite{GabdullinMakarenkoShevtsova2018}~--- directly to all the functions in~$\mathcal G$ (see~\eqref{KatzPetrovGproperties}).  On the other hand, Bikelis' inequality~\eqref{BikelisIneqMin} trivially follows from~\eqref{PetrovIneq}  with 
$$
g(u)=g_*(u;(1+|x|)B_n)=\min\Big\{1,\frac{|u|}{(1+|x|)B_n}\Big\}\in\mathcal G.
$$
Moreover, inequality~\eqref{PetrovIneq} with $g(u)=|u|^\d\in\mathcal G$ also yields Nagaev--Bikelis inequality~\eqref{NagaevBikelisIneq} with $K_0\le4A$. However, the numerical optimization of the constant $A$ with the concrete function $g(u)=|u|^\d$ (which is not an extremal in~\eqref{PetrovIneq}), in fact, leads to sharper upper bounds for $A$ which coincides with $K_0$ in this case (see table~\ref{TabBEnonuni1+x_struct}) than those that can be obtained for the universal constant $A$ (with the extremal function $g=g_*$).

Let us also mention that, in 2001, Chen and Shao~\cite{ChenShao2001} reproved Bikelis' inequality~\eqref{BikelisIneq}  by Stein's method; moreover, the authors of~\cite{ChenShao2001} refer to Bikelis' work~\cite{Bikelis1966}, but cite only weaker inequality~\eqref{NagaevBikelisIneq} stating erroneously that results of~\cite{Bikelis1966} are of a less general character and hold true only under the assumption of finiteness of third-order moments of random summands.

The value of the constant  $A$ also remained unknown for a long time. It's first upper bounds were obtained only in 2005--2007 by Neammanee and Thongtha~\cite{Neammanee2005,ThongthaNeammanee2007,NeammaneeThongtha2007}. Then they were improved by Korolev and Popov~\cite{KorolevPopov2012} to the presently best known bounds:  $A\le 39.32$ in the i.i.d. case and $A\le 47.65$ in the general situation. Moreover, in~\cite{KorolevPopov2012} it is also shown that the following improved estimates hold for large values of the argument $|x|\ge10$: $A\le 24.13$ in the i.i.d. case and $A\le 29.62$ in the general situation.

Let us also note that, in 2011, Gavrilenko, Nefedova, and Shevtsova~\cite{Gavrilenko2011,NefedovaShevtsova2011DAN,Shevtsova2013Inf,Shevtsova2017book} suggested structural improvements of Nagaev--Bikelis inequality~\eqref{NagaevBikelisIneq} in the following form:
\begin{equation}\label{BEnonuni1+x_struct}
\sup_{x\in\R}(1+|x|^{2+\d})\Delta_n(x)\le \inf_{s\ge0}K_s(\d)\big(\lyapd+sT_{2+\d,n}\big),
\end{equation}
where
$$
T_{2+\d,n}:= \frac1{B_n^{2+\d}}\sum_{j=1}^n\sigma_j^{2+\d}\ \le\ \frac1{B_n^{2+\d}}\sum_{j=1}^n\E|X_j|^{2+\d} =\lyapd
$$
and, of course, $K_s(\d)\le K_0(\d)$ for all $s\ge0$. However, values of the constants $K_s(\d)$ for $s>0$ turn to be strictly less than for $s=0$, which makes estimate~\eqref{BEnonuni1+x_struct} more favorable than the classical Nagaev--Bikelis inequality~\eqref{NagaevBikelisIneq} for large values of the ratio $\lyapd/T_{2+\d,n}$ (which is never less than one and may be infinitely large). The best known upper bounds for the constants  $K_s(\d)$ are obtained in~\cite{Shevtsova2013Inf,Shevtsova2017book} and are cited in table~\ref{TabBEnonuni1+x_struct} for some $s\in[0,1]$ and $\d\in(0,1]$, where $s_1(\d)$ is the optimal value of $s\ge0$, that minimizes $K_s(\d)$ (within the method used), so that  $K_s(\d)=K_{s_1}(\d)$ for $s>s_1(\d)$.

As regards lower bounds for the constants $K_0(\d)$ and $A$, presently a couple of lower bounds is known only for $K_0(1)$. The first one follows from Chistyakov's result~\cite{Chistyakov1990}:
$$
K_0(1)\ge \lim_{|x|\to\infty}\limsup_{\ell\to0}\sup_{n,X_1,\ldots,X_n\colon L_{3,n}=\ell}|x|^3\Delta_n(x)/\ell =1.
$$
The second one was obtained in a recent paper~\cite{Pinelis2013}  by Pinelis who considered $n=1$, $\Prob(X_1=1-p)=p=1-\Prob(X_1=-p)$, ${x=1-p}$ and $p=0.08$ and improved the above bound to
$$
K_0(1)\ge \sup_{x\in\R,X_1}
\left (1+|x|^3\right )\Delta_1(x)\frac{\sigma_1^3}{\E|X_1|^3}>1.0135\ldots
$$

\section{Main results}

Using Pinelis' lower bound for $K_0(1)$ and inequality~\eqref{K0<=A*2^(1+d)} it is easy to obtain a lower bound for the constant $A$ in Bikelis'~\eqref{BikelisIneqInt}, \eqref{BikelisIneqMin}, \eqref{BikelisIneq} and Petrov's~\eqref{PetrovIneq} inequalities in the following form:
$$
A\ge\sup_{\d\in(0,1]}K_0(\d)/2^{1+\d}\ge K_0(1)/4>0.2533.
$$
However, one can act more delicate (similarly to Pinelis~\cite{Pinelis2013}) and obtain a sharper lower bound. 

\begin{theorem}\label{ThBikelisPetrovConstLowBound}
For the absolute constant $A$ in Bikelis'~\eqref{BikelisIneqInt}, \eqref{BikelisIneqMin}, \eqref{BikelisIneq} and Petrov's~\eqref{PetrovIneq} inequalities, also in the i.i.d. case, we have
$$
A \ge \big(1+\sqrt{p^{-1}-1}\,\big)^2 \big(\Phi(\sqrt{p^{-1}-1}\,)-1+p\big)\Big|_{p=0.15}>1.6153.
$$
\end{theorem}


Similar reasoning leads to the following lower bounds for the constants $K_0(\d)$ in Nagaev--Bikelis inequality~\eqref{NagaevBikelisIneq} and $K_s(\d)$ in~\eqref{BEnonuni1+x_struct} with arbitrary $\d\in[0,1]$ (observe that inequalities~\eqref{NagaevBikelisIneq}, \eqref{BEnonuni1+x_struct} hold true also for $\d=0$ with $\lyapd=T_{2+\d,n}=1$, as it follows, say, from~\eqref{PetrovIneq} with $g(u)\equiv1$).

\begin{theorem}\label{ThNagaevBikelisConstLowBound}
For the constants $K_0(\d)$ in~\eqref{NagaevBikelisIneq} and $K_s(\d)$ in~\eqref{BEnonuni1+x_struct} for every $s\ge0$ and $\d\in[0,1]$, also in the i.i.d. case, we have
\begin{equation}\label{NBstructConst(delta)>=}
K_s(\d)\ge\sup_{0<p<1,\,q=1-p}  q^{\d/2} \cdot \frac{p^{1+\d/2}+q^{1+\d/2}}{p^{1+\d}+q^{1+\d}+s(pq)^{\d/2}} \abs{1-\frac1p\Phi\left(-\sqrt{\frac{q}p}\,\right )},
\end{equation}
in particular,
\begin{equation}\label{NBconst(delta)>=}
K_0(\d)\ge\sup_{0<p<1,\,q=1-p}  q^{\d/2} \cdot \frac{p^{1+\d/2}+q^{1+\d/2}}{p^{1+\d}+q^{1+\d}} \abs{1-\frac1p\Phi\left(-\sqrt{\frac{q}p}\,\right )},\quad \d\in[0,1],
\end{equation}
\begin{equation}\label{NBstructConst(delta)>=uniform(p->0)}
K_s(\d)\ge
\begin{cases}
1/(1+s),&\d=0,\\
1,&\d\in(0,1],
\end{cases}
\quad s\ge0.
\end{equation}
\end{theorem}

The lower bound  in~\eqref{NBstructConst(delta)>=uniform(p->0)} is obtained by letting $p\to0+$ in \eqref{NBstructConst(delta)>=}.  It turns out that $p\to0+$ is indeed an extremal for either $\d=0$ or $s>0$ (the numerically optimal values of $p$ are very close to zero), so we leave lower bounds in~\eqref{NBstructConst(delta)>=uniform(p->0)} as finite ones for $\d=0$ and all $s\ge0$ or $s>0$ and all $\d\in[0,1]$, while an accurate optimization in~\eqref{NBconst(delta)>=}  with respect to $p\in(0,1)$ for fixed  $\d\in(0,1]$ leads to sharper lower bounds for $K_0(\d)$ given in table~\ref{TabNBconst(delta)>=} in the second row. The values of the minorant~\eqref{NBconst(delta)>=}  in table~\ref{TabNBconst(delta)>=} are accompanied with the corresponding values of $p$ (in the third row) close to the extremal ones which guarantee the announced lower bounds.

\begin{table}[h]
\centering
\begin{tabular}{||c||c|c|c|c|c|c|c|c|c|c|c||}
\hline
$\d$&0&0.1&0.2&0.3&0.4&0.5&0.6&0.7&0.8&0.9&1
\\\hline
$K_0\ge$ &1&1.0061&1.0108&1.0139&1.0158&1.0167&1.0168&1.0164&1.0157&1.0147&1.0135
\\\hline
$p$&0&0.06&0.066&0.07&0.074&0.076&0.08&0.08&0.08&0.08&0.08
\\\hline
\end{tabular}
\caption{Lower bounds for the constants $K_0(\delta)$ from~\eqref{NagaevBikelisIneq}, constructed with respect to formula~\eqref{NBconst(delta)>=}, for some $\delta\in[0,1]$.}
\label{TabNBconst(delta)>=}
\end{table}

\section{Proofs}

\subsection{Proof of Theorem~\ref{ThBikelisPetrovConstLowBound}}
To construct a lower bound for the constant $A$ consider Petrov's inequality~\eqref{PetrovIneq}  with
$$
g(u)\equiv1\in\mathcal G\quad\text{and}\quad n=1.
$$
Then we have
$$
A\ge\sup_{x\in\R,X_1}(1+|x|)^2\Delta_1(x),
$$
where the least upper bound is taken with respect to $x\in\R$ and all distributions of the r.v. $X_1$ with $\E X_1=0$, $\E X_1^2\in(0,\infty)$. Now letting
$$
\Prob\big(X_1=\sqrt{q/p}\,\big)=p=1-\Prob\big(X_1=-\sqrt{p/q}\,\big),\quad x=\sqrt{q/p},\quad q=1-p,\quad p\in(0,1),
$$
we obtain
$$
A\ge \sup_{p\in(0,1),\,q=1-p}\big(1+\sqrt{q/p}\,\big)^2 \abs{q-\Phi(\sqrt{q/p})}.
$$
The announced lower bound follows by taking here $p=0.15$.

\subsection{Proof of Theorem~\ref{ThNagaevBikelisConstLowBound}}
To construct lower bounds for the constants $K_s(\d)$ consider inequality~\eqref{BEnonuni1+x_struct} with
$$
n=1,\quad\Prob\big(X_1=\sqrt{q/p}\,\big)=p=1-\Prob\big(X_1=-\sqrt{p/q}\,\big),\quad x=\sqrt{q/p},\quad q=1-p,\quad p\in(0,1).
$$
Then we have
$$
\E X_1=0,\quad \sigma_1^2=1,\quad \E|X_1|^{2+\d}= \frac{p^{1+\d}+q^{1+\d}}{(pq)^{\d/2}},
$$
and
\begin{multline*}
K_s(\d)\ge \sup_{x\in\R,X_1}
\left (1+|x|^{2+\d}\right) \frac{\Delta_1(x)\sigma_1^{2+\d}}{\E|X_1|^{2+\d}+s\sigma_1^{2+\d}} 
\\
\ge\sup_{0<p<1,\,q=1-p} \left(1+\Big(\frac{q}p\Big)^{1+\d/2}\right) \abs{q-\Phi\left (\sqrt{\frac qp}\,\right)} \frac{(pq)^{\d/2}}{p^{1+\d}+q^{1+\d}+s(pq)^{\d/2}}
\\
=\sup_{0<p<1,\,q=1-p}  q^{\d/2} \cdot \frac{p^{1+\d/2}+q^{1+\d/2}}{p^{1+\d}+q^{1+\d}+s(pq)^{\d/2}} \abs{1-\frac1p\Phi\left(-\sqrt{\frac{q}p}\,\right )}
\end{multline*}
for all $\d\in[0,1]$ and $s\ge0$. Now letting $p\to0+$ and taking into account that 
$$
\Phi(-x)\le \frac1{\sqrt{2\pi}x}e^{-x^2/2}, \quad x>0,
$$
we obtain a lower bound
$$
K_s(\d)\ge \lim_{p\to0}\abs{1-\frac1p\Phi\left(-\sqrt{\frac{1-p}p}\,\right )} =  1-\lim_{x\to\infty}(1+x^2)\Phi(-x)=1,
$$
universal for all $s\ge0$ and $\d\in[0,1]$, while an accurate optimization with respect to $p\in(0,1)$ for fixed $\d\in(0,1]$ and $s=0$ leads to sharper lower bounds for $K_0(\d)$ given in table~\ref{TabNBconst(delta)>=}.

\AcknowledgementSection
This work was supported by the Russian Foundation for Basic Research (project No.\,19-07-01220-a) and by the Ministry for Education and Science of Russia (grant No.\,MD--189.2019.1).

\end{document}